\documentclass[11pt]{article}
\usepackage{graphicx}
\usepackage{amsfonts}
\usepackage{nath}
  \def\fH{{\cal H}}

\usepackage{theorem}
\newtheorem{Lem}{Lemma}[section]
\newtheorem{Def}[Lem]{Definition}
\newtheorem{The}[Lem]{Theorem}
\newtheorem{Prop}[Lem]{Proposition}

\newtheorem{Rem}[Lem]{Remark}

\newcommand{\qed}{\hbox{\rule{6pt}{6pt}}}

\begin{document}
\title{Unitarily invariant norm inequalities for some means}
\author{Shigeru Furuichi\footnote{E-mail:furuichi@chs.nihon-u.ac.jp} \\
$^1${\small Department of Information Science,}\\
{\small College of Humanities and Sciences, Nihon University,}\\
{\small Sakurajyousui, Setagaya-ku, Tokyo, 156-8550, Japan}}
\date{}
\maketitle
{\bf Abstract.} 
We introduce some symmetric homogeneous means, and then show unitarily invariant norm inequalities for them, applying the method established by Hiai and Kosaki. Our new inequalities give the tighter bounds of the logarithmic mean than the inequalities given by Hiai and Kosaki. Some properties and norm continuities in parameter for our means are also discussed.

\vspace{3mm}


\vspace{3mm}
\vspace{3mm}

%%%%%%%%%%%%%%%%%%%%%%%%%%%%%%%%%%%%%%%%%%%%%%%%%%%%%%%%%%%%%%%%%%%%%%%%%%%%%%%%%%%%%%%%%%%%%%%%%%%%%%%%%%%%%%%%%%%%%%%%%%%%%%%%%%%%%%%%%%%%%%%%%%%%%%%%%%%%%%%
%%%%%%%%%%%%%%%%%%%%%%%%%%%%%%%%%%%%%%%%%%%%%%%%%%%%%%%%%%%%%%%%%%%%%%%%%%%%%%%%%%%%%%%%%%%%%%%%%%%%%%%%%%%%%%%%%%%%%%%%%%%%%%%%%%%%%%%%%%%%%%%%%%%%%%%%%%%%%%%
%%%%%%%%%%%%%%%%%%%%%%%%%%%%%%%%%%%%%%%%%%%%%%%%%%%%%%%%%%%%%%%%%%%%%%%%%%%%%%%%%%%%%%%%%%%%%%%%%%%%%%%%%%%%%%%%%%%%%%%%%%%%%%%%%%%%%%%%%%%%%%%%%%%%%%%%%%%%%%%
%%%%%%%%%%%%%%%%%%%%%%%%%%%%%%%%%%%%%%%%%%%  Section1  %%%%%%%%%%%%%%%%%%%%%%%%%%%%%%%%%%%%%%%%%%%%%%%%%%%%%%%%%%%%%%%%%%%%%%%%%%%%%%%%%%%%%%%%%%%%%%%%%%%%%%%%
%%%%%%%%%%%%%%%%%%%%%%%%%%%%%%%%%%%%%%%%%%%%%%%%%%%%%%%%%%%%%%%%%%%%%%%%%%%%%%%%%%%%%%%%%%%%%%%%%%%%%%%%%%%%%%%%%%%%%%%%%%%%%%%%%%%%%%%%%%%%%%%%%%%%%%%%%%%%%%%
%%%%%%%%%%%%%%%%%%%%%%%%%%%%%%%%%%%%%%%%%%%%%%%%%%%%%%%%%%%%%%%%%%%%%%%%%%%%%%%%%%%%%%%%%%%%%%%%%%%%%%%%%%%%%%%%%%%%%%%%%%%%%%%%%%%%%%%%%%%%%%%%%%%%%%%%%%%%%%%
%%%%%%%%%%%%%%%%%%%%%%%%%%%%%%%%%%%%%%%%%%%%%%%%%%%%%%%%%%%%%%%%%%%%%%%%%%%%%%%%%%%%%%%%%%%%%%%%%%%%%%%%%%%%%%%%%%%%%%%%%%%%%%%%%%%%%%%%%%%%%%%%%%%%%%%%%%%%%%%
\section{Introduction}
In the previous paper, we derived the tight bounds for logarithmic mean in the case of Frobenius norm, inspired by the work of Zou in \cite{Zou}.  
\begin{The} {\bf (\cite{FY})} \label{FY_theorem2013}
For any matrices $S,T,X$ with $S,T\geq 0$, $m_1 \geq 1$, $m_2 \geq 2$ and Frobenius norm $\lVert \cdot \rVert_F$, the following inequalities hold.
\begin{eqnarray*}
\hspace*{-12mm }\frac{1}{m_1}\lVert\sum_{k=1}^{m_1}S^{k/(m_1+1)}XT^{(m_1+1-k)/(m_1+1)}\rVert_F &\leq&\frac{1}{m_1} \lVert \sum_{k=1}^{m_1} S^{(2k-1)/2m_1} XT^{(2m_1-(2k-1))/2m_1} \rVert_F\\
&\leq& \lVert \int_0^1 S^{\nu}XT^{1-\nu}d\nu \rVert_F \\
&\leq& \frac{1}{m_2} \lVert\sum_{k=0}^{m_2} S^{k/m_2}XT^{(m_2-k)/m_2}-\frac{1}{2}(SX+XT) \rVert_F \\
&\leq& \frac{1}{m_2} \lVert \sum_{k=0}^{m_2-1} S^{k/(m_2-1)}XT^{(m_2-1-k)/(m_2-1)}\rVert_F.
\end{eqnarray*}
\end{The}
Our bounds for the logarithmic mean have improved the famous results by Hiai and Kosaki \cite{HK1,HK2} in the special case, since Frobenius norm is one of unitarily invariant norms.
\begin{The} {\bf (\cite{HK1,HK2})}   \label{HK_theorem1999}
For any bounded linear operators $S,T,X$ with $S,T \geq 0$, $m_1 \geq 1$, $m_2 \geq 2$ and any  unitarily invariant norm $\ltriple| \cdot \rtriple|$, the following inequalities hold.
\begin{eqnarray*}
\ltriple| S^{1/2}XT^{1/2}\rtriple| &\leq& \frac{1}{m_1}\ltriple|\sum_{k=1}^{m_1}S^{k/(m_1+1)}XT^{(m_1+1-k)/(m_1+1)}\rtriple| \leq \ltriple| \int_0^1 S^{\nu}XT^{1-\nu}d\nu \rtriple|\\
&\leq& \frac{1}{m_2} \ltriple| \sum_{k=0}^{m_2-1} S^{k/(m_2-1)}XT^{(m_2-1-k)/(m_2-1)}\rtriple| \leq \frac{1}{2} \ltriple| SX +XT\rtriple|.
\end{eqnarray*}
\end{The}
In this paper, we give the tighter bounds for the logarithmic mean than those by Hiai and Kosaki \cite{HK1,HK2} for every unitarily invariant norm. That is, we give the generalized results of Theorem \ref{FY_theorem2013} for the unitarily invariant norm. For this purpose, we firstly introduce two quantities.
%non-symmetric homogeneous means (See Appendix A.1 of \cite{Hiai-Kosaki-book}).
\begin{Def} \label{def-P-Q}
For $\alpha \in \mathbb{R}$ and $x, y >0$, we set
\[ P_{\alpha} ( x,y ) \equiv 
\left\{ \begin{array}{l}
\frac{\alpha x^{\alpha} ( x - y )}{x^{\alpha}  - y^{\alpha} },\,\,\,\,( x \neq y )\\
\,\,\,\,\,\,x,\,\,\,\,\,\,\,\,\,\,\,\,\,\,\,\,\,\,\,\,\,\,\,\,\, ( x = y )
\end{array} \right.
\,\,\,and\,\,\,\, Q_{\alpha} ( x,y ) \equiv 
\left\{ \begin{array}{l}
\frac{\alpha y^{\alpha} ( x - y )}{x^{\alpha} - y^{\alpha}},\,\,\,\,( x \neq y )\\
\,\,\,\,\,\,x,\,\,\,\,\,\,\,\,\,\,\,\,\,\,\,\,\,\,\,\,\,\,\,\,\,( x = y ).
\end{array} \right.
\]
\end{Def}

We note that we have the following bounds of logarithmic mean with the above two means (See Appendix in the paper \cite{FY}):
\[\left\{ \begin{array}{l}
Q_{1/m}(x,y) < LM( x,y ) < P_{1/m}( x,y ),( if\,\,\,x > y ),\\
P_{1/m}( x,y ) < LM( x,y ) < Q_{1/m}( x,y ),( if\,\,\,x < y ),
\end{array} \right.\,\,\,\]
where the logarithmic mean is defined by
\begin{equation}
LM( x,y ) \equiv \left\{ \begin{array}{l}
\frac{x - y}{\log x - \log y},\,\,\,\,( x \ne y )\\
\,\,\,\,x,\,\,\,\,\,\,\,\,\,\,\,\,\,\,\,\,\,\,\,\,\,\,\,\,\,\,\,\,\,\,( x = y).
\end{array} \right.\,
\end{equation}
We here define a few symmetric homogeneous means using $P_{\alpha}(x,y)$ and $Q_{\alpha}(x,y)$ in the following way.
\begin{Def} \label{def-A-L-G-H}
\begin{itemize}
\item[(i)] For $\vert \alpha\vert \leq 1$ and $x \neq y$, we define, 
\[
A_{\alpha}( x,y ) \equiv \frac{1}{2} P_{\alpha}( x,y) + \frac{1}{2}Q_{\alpha }( x,y ) = \frac{\alpha ( x^{\alpha } + y^{\alpha } )( x - y )}{2( x^{\alpha } - y^{\alpha })},
\]
\item[(ii)] For $\alpha \in \mathbb{R}$ and $x \neq y$, we define, 
\[
L_{\alpha }( x,y) \equiv \frac{P_{\alpha }( x,y ) - Q_{\alpha }( x,y )}{\log P_{\alpha }( x,y ) - \log Q_{\alpha }( x,y )} = LM( x,y ),
\]
\item[(iii)] For $\vert \alpha\vert \leq 2$ and $x \neq y$, we define, 
\[
G_{\alpha }( x,y ) \equiv \sqrt {P_{\alpha }( x,y )Q_{\alpha }( x,y )}  = \frac{\alpha ( xy )^{\alpha /2}( x - y )}{x^{\alpha } - y^{\alpha }},
\]
\item[(iv)] For $\vert \alpha\vert \leq 1$ and $x \neq y$, we define, 
\[
H_{\alpha }( x,y ) \equiv \frac{ 2 P_{\alpha }( x,y ) Q_{\alpha }( x,y )}{P_{\alpha }( x,y ) + Q_{\alpha }( x,y )} = \frac{2\alpha ( xy )^{\alpha }}{x^{\alpha } + y^{\alpha }}\frac{( x - y )}{x^{\alpha } - y^{\alpha }},
\]
\end{itemize}
and we also set  
${A_\alpha }\left( {x,y} \right) = {L_\alpha }\left( {x,y} \right) = {G_\alpha }\left( {x,y} \right) = {H_\alpha }\left( {x,y} \right)=x$ for $x=y$.
\end{Def}

We have the following relations for the above means:
\begin{eqnarray*}
&&A_1(x,y) = AM(x,y) \equiv \frac{1}{2}(x+y),\quad A_0(x,y) =
  \lim_{\alpha \rightarrow 0} A_{\alpha}(x,y) = LM(x,y),\\
&&G_{0}(x,y) = \lim_{\alpha \rightarrow 0} G_{\alpha}(x,y) =LM(x,y), \quad G_1(x,y) = GM(x,y) \equiv \sqrt{xy}, \\
&&G_2(x,y) = HM(x,y) \equiv \frac{2xy}{x+y},\quad H_0(x,y) = \lim_{\alpha \rightarrow 0}H_{\alpha}(x,y) =LM(x,y),\\
&& H_{1/2}(x,y) = GM(x,y), \quad H_1(x,y)=HM(x,y)
\end{eqnarray*} 
and $H_{\alpha}(x,y) =G_{2\alpha}(x,y)$.
In addition, the above means are written as the following geometric
bridges:
\begin{eqnarray*}
&&A_{\alpha}(x,y) = [B_{\alpha}(x,y)]^{\alpha} [S_{\alpha}(x,y)]^{1-\alpha}, \quad
L_{\alpha}(x,y) = [E_{\alpha}(x,y)]^{\alpha} [S_{\alpha}(x,y)]^{1-\alpha}, \\
&&G_{\alpha}(x,y) = [GM(x,y)]^{\alpha} [S_{\alpha}(x,y)]^{1-\alpha}, \quad
H_{\alpha}(x,y) = [D_{\alpha}(x,y)]^{\alpha} [S_{\alpha}(x,y)]^{1-\alpha}, 
\end{eqnarray*}
where
$$
S_{\alpha}(x,y) \equiv \left( \frac{\alpha(x-y)}{x^{\alpha}-y^{\alpha}}\right)^{\frac{1}{1-\alpha}},\quad B_{\alpha}(x,y) \equiv \left(\frac{x^{\alpha}+y^{\alpha}}{2} \right)^{\frac{1}{\alpha}}
$$
and
$$
D_{\alpha}(x,y) \equiv\left(\frac{2x^{\alpha} y^{\alpha}}{x^{\alpha} +y^{\alpha}} \right)^{\frac{1}{\alpha}},  \quad E_{\alpha}(x,y) \equiv\left( \frac{x^{\alpha}-y^{\alpha}}{\alpha(\log x -\log y)}\right)^{\frac{1}{\alpha}}. 
$$
$S_{\alpha}(x,y) $ and  $B_{\alpha}(x,y)$ are called  Stolarsky mean and binomial mean, respectively.

In the previous paper \cite{FY}, as tight bounds of logarithmic mean, the scalar inequalities were shown
$$
G_{1/m}(x,y) \leq LM(x,y), \quad (m \geq 1), \qquad LM(x,y) \leq A_{1/m}(x,y), \quad (m \geq 2)
$$
which equivalently implied Frobenius norm inequalities (Theorem \ref{FY_theorem2013}). See Theorem 2.2 and Theorem 3.2 in \cite{FY} for details. In this paper, we give unitarily invariant norm inequalities which are general results including Frobenius norm inequalities as a special case.

%%%%%%%%%%%%%%%%%%%%%%%%%%%%%%%%%%%%%%%%%%%%%%%%%%%%%%%%%%%%%%%%%%%%%%%%%%%%%%%%%%%%%%%%%%%%%%%%%%%%%%%%%%%%%%%%%%%%%%%%%%%%%%%%%%%%%%%%%%%%%%%%%%%%%%%%
%%%%%%%%%%%%%%%%%%%%%%%%%%%%%%%%%%%%%%%%%%%%%%%%%%%%%%%%%%%%%%%%%%%%%%%%%%%%%%%%%%%%%%%%%%%%%%%%%%%%%%%%%%%%%%%%%%%%%%%%%%%%%%%%%%%%%%%%%%%%%%%%%%%%%%%%
%%%%%%%%%%%%%%%%%%%%%%%%%%%%%%%%%%%%%%%%%%%%%%%%%%%%%%%%%%%%%%%%%%%%%%%%%%%%%%%%%%%%%%%%%%%%%%%%%%%%%%%%%%%%%%%%%%%%%%%%%%%%%%%%%%%%%%%%%%%%%%%%%%%%%%%%%
%%%%%%%%%%%%%%%%%%%%%%%%%%%%%%%%%%%%%%%%%%%%%%%%%%%%%%%%%%%%%%%%%%%%%%%%%%%%%%%%%%%%%%%%%%%%%%%%%%%%%%%%%%%%%%%%%%%%%%%%%%%%%%%%%%%%%%%%%%%%%%%%%%%%%%%%
%%%%%%%%%%%%%%%%%%%%%%%%%%%%%%%%%%%%%%%%%%%%%%%%%%%%%%%%%%%%%%%%%%%%%%%%%%%%%%%%%%%%%%%%%%%%%%%%%%%%%%%%%%%%%%%%%%%%%%%%%%%%%%%%%%%%%%%%%%%%%%%%%%%%%%%%
%%%%%%%%%%%%%%%%%%%%%%%%%%%%%%%%%%%%%%%%%%%%%%%%%%%%%%%%%%%%%%%%%%%%%%%%%%%%%%%%%%%%%%%%%%%%%%%%%%%%%%%%%%%%%%%%%%%%%%%%%%%%%%%%%%%%%%%%%%%%%%%%%%%%%%%%
%%%%%%%%%%%%%%%%%%%%%%%%%%%%%%%%%%%%%%%%%%%%%%%%%%%%%%%%%%%%%%%%%%%%%%%%%%%%%%%%%%%%%%%%%%%%%%%%%%%%%%%%%%%%%%%%%%%%%%%%%%%%%%%%%%%%%%%%%%%%%%%%%%%%%%%%%
%%%%%%%%%%%%%%%%%%%%%%%%%%%%%%%%%%%%%%%%%%%%%%%%%%%%%%%%%%%%%%%%%%%%%%%%%%%%%%%%%%%%%%%%%%%%%%%%%%%%%%%%%%%%%%%%%%%%%%%%%%%%%%%%%%%%%%%%%%%%%%%%%%%%%%%%

\section{Unitarily invariant norm inequalities}
To obtain unitarily invariant norm inequalities, we apply the method established by Hiai and Kosaki \cite{HK2,Hiai-Kosaki-book,Hiai2010,Kosaki}. 

\begin{Def} \label{Def-SHM-section2}
A continuous positive real function $M(x,y)$ for $x, y >0$ is called a symmetric homogeneous mean if the function $M$ satisfies the following properties:
\begin{itemize}
\item[(a)] $M(x,y) = M(y,x)$.
\item[(b)] $M(cx,cy)=cM(x,y)$ for $c>0$.
\item[(c)] $M(x,y)$ is non-decreasing in $x, y$.
\item[(d)] $\min\{x, y\} \leq M(x, y) \leq \max\{x, y\}$.
\end{itemize}
\end{Def}
The functions $A_{\alpha}(x,y), L_{\alpha}(x,y), G_{\alpha}(x,y), H_{\alpha}(x,y)$ defined in Definition \ref{def-A-L-G-H} are symmetric homogeneous means.
%(See Appendix.)
We give powerful theorem to obtain unitarily invariant norm inequalities.
In the references \cite{HK2,Hiai-Kosaki-book,Hiai2010,Kosaki}, another equivalent conditions were given. However here we give minimum conditions to obtain our results in this paper. Throughout this paper, we use the symbol $B(\fH)$ as the set of all bounded linear operators on a separable Hilbert space $\fH$. We also use the notation $K \geq 0$ if $K\in B(\fH)$ satisfies $\langle Kx,x\rangle \geq 0$ for all $x \in \fH$ (then $K$ is called a positive operator).

\begin{The} {\bf (\cite{HK2,Hiai-Kosaki-book,Hiai2010,Kosaki})} \label{Hiai-Kosaki-theorem}
For two symmetric homogeneous means $M$ and $N$, the following conditions are equivalent:
\begin{itemize}
\item[(i)] $\ltriple| M(S,T)X \rtriple| \leq \ltriple| N(S,T)X \rtriple|$ for any $S, T, X \in B(\fH)$ with $S, T \geq 0$ and for any unitarily invariant norm $\ltriple| \cdot \rtriple|$.
\item[(ii)] The function $M(e^t,1) / N(e^t,1)$ is positive definite function on $\mathbb{R}$ (then we denote $M \preceq N$), where the positive definiteness of a real continuous function $\phi$ on $\mathbb{R}$ means that $[\phi(t_i-t_j)]_{i,j=1,\cdots,n}$ is positive definite for any $t_1,\cdots,t_n \in \mathbb{R}$ with any $n \in \mathbb{N}$. 
\end{itemize}
\end{The}

Thanks to Theorem \ref{Hiai-Kosaki-theorem}, our task to obtain unitarily invariant norm inequalities in this paper is to show the relation $M \preceq N$ which is stronger than the usual scalar inequalities $M \leq N$.
That is, $M(s,t) \preceq N(s,t)$ implies $M(s,t) \leq N(s,t)$. 

We firstly give monotonicity of three means $H_{\alpha}(x,y)$, $G_{\alpha}(x,y)$ and
$A_{\alpha}(x,y)$ for the parameter $\alpha \in \mathbb{R}$.
Since we have $H_{-\alpha}(x,y) =H_{\alpha}(x,y)$, $G_{-\alpha}(x,y) =G_{\alpha}(x,y)$ and $A_{-\alpha}(x,y) =A_{\alpha}(x,y)$, we consider the case $\alpha \geq 0$.
Then we have the following proposition.
\begin{Prop} \label{prop-monotonicity-section2}
\begin{itemize}
\item[(i)] If $0 \leq \alpha < \beta \leq 1$, then $H_{\beta} \preceq H_{\alpha}$.
\item[(ii)] If $0 \leq \alpha < \beta \leq 2$, then $G_{\beta} \preceq G_{\alpha}$.
\item[(iii)] If $0 \leq \alpha < \beta \leq 1$, then $A_{\alpha} \preceq A_{\beta}$.
\end{itemize}
\end{Prop}

{\it Proof:}
\begin{itemize}
\item[(i)] We calculate
\[
\frac{H_{\beta }( e^t,1 )}{H_{\alpha }( e^t,1 )} 
= \frac{ 2\beta  e^{\beta t} ( e^t - 1 )}{e^{2\beta t} - 1} 
\cdot \frac{e^{2\alpha t} - 1}{2\alpha e^{\alpha t}( e^t - 1 )} 
= \frac{\beta }{\alpha } \cdot \frac{e^{\beta t}( e^{2\alpha t} - 1 )}{e^{\alpha t}( e^{2\beta t} - 1 )} 
= \frac{\beta }{\alpha }\frac{ \sinh \alpha t}{\sinh \beta t}.
\]
This is a positive definite function for the case $\alpha < \beta$, so that
we have $H_{\beta} \preceq H_{\alpha}$.
\item[(ii)] The similar calculation
\[
\frac{G_{\beta }( e^{2t},1 )}{G_{\alpha }( e^{2t},1 )} 
= \frac{ 2\beta e^{\beta t}( e^{2t} - 1 )}{e^{2\beta t}- 1} \cdot \frac{e^{2\alpha t} - 1}{2\alpha e^{\alpha t}( e^{2t} - 1 )} 
= \frac{\beta }{\alpha } \cdot \frac{e^{\beta t} (e^{2\alpha t} - 1 )}{e^{\alpha t}( e^{2\beta t} - 1 )} 
= \frac{\beta }{\alpha } \cdot\frac{{\sinh \alpha t}}{{\sinh \beta t}}
\]
implies $G_{\beta} \preceq G_{\alpha}$.
\item[(iii)] Since we have
\[
\frac{A_{\alpha }( e^{2t},1 )}{A_{\beta }( e^{2t},1)} 
= \frac{\alpha }{\beta } \cdot \frac{{\sinh \beta t\cosh \alpha t}}{{\cosh \beta t\sinh \alpha t}},\]
we calculate by the formula $\sinh(x) = 2\cosh\left(\frac{x}{2}\right)\sinh\left(\frac{x}{2}\right)$ repeatedly
\begin{eqnarray*}
\frac{\sinh \beta t\cosh \alpha t}{\cosh \beta t\sinh \alpha t} - 1 &=& \frac{\sinh ( \beta  - \alpha  )t}{\cosh \beta t\sinh \alpha t}
=\frac{2\cosh (\frac{\beta -\alpha}{2}t )\sinh (\frac{\beta -\alpha}{2}t)}{\cosh \beta t\sinh \alpha t}\\
&=&\lim_{n\rightarrow \infty}\frac{2^n \prod_{k=1}^n\cos (\frac{\beta -\alpha}{2^k}t )\sinh (\frac{\beta -\alpha}{2^n}t ) }{\cosh \beta t\sinh \alpha t}.
\end{eqnarray*}
From Proposition 4 in \cite{BK}, the sufficient condition that the function
$\frac{\prod_{k=1}^n \cosh\left(\frac{\beta-\alpha}{2^k}t\right)}{\cosh\beta t}$ is positive definite, is $\sum_{k=1}^n\frac{\beta -\alpha}{\beta 2^k} \leq 1$, i.e., $(\beta-\alpha)(1-2^{-n}) \leq \beta$.
The sufficient condition that the function $\frac{\sinh\left(\frac{\beta-\alpha}{2^n}t\right)}{\sinh\alpha t}$ is positive definite, is $\frac{\beta-\alpha}{2^n}\leq \alpha$. When $n \rightarrow \infty$, both conditions become to $0 \leq \alpha$, which satisfies the assumption of this proposition. Thus we conclude $A_{\alpha}  \preceq A_{\beta}$. 

%%%%%%% The below is the proof by Fumio Hiai%%%%%%%
%If $\alpha < \beta \leq 2\alpha$, then $0 < \beta -\alpha \leq \alpha$ so that
%$\frac{{\sinh \left( {\beta  - \alpha } \right)t}}{{\sinh \alpha t}}$ is a positive %definite function. The function $\frac{1}{{\cosh \beta t}}$ is also positive definite.
%Therefore we have $A_{\alpha} \preceq A_{\beta}$ for the case $\alpha < \beta \leq 2 \alpha$. For the general case $\alpha < \beta \leq 1$, we can show
%$A_{\alpha} \preceq A_{\beta}$ by the similar argument in Theorem 2.1 of the paper \cite{HK2}.

\end{itemize}

\hfill \qed

It may be notable that (iii) of the above proposition can be proven by the similar argument in Theorem 2.1 of the paper \cite{HK2}.

Next we give the relation among four means $H_{\alpha}(x,y)$,  $G_{\alpha}(x,y)$,  $L_{\alpha}(x,y)$, and  $A_{\alpha}(x,y)$.
\begin{Prop} \label{prop-ineq-section2}
For any $S, T, X \in B(\fH)$ with $S, T \geq 0$, $\vert \alpha  \vert \leq 1$ and any unitarily invariant norm $\ltriple| \cdot \rtriple|$, we have
$$
\ltriple| H_{\alpha}(S, T)X \rtriple| \leq 
\ltriple| G_{\alpha}(S, T)X \rtriple| \leq 
\ltriple| L_{\alpha}(S, T)X \rtriple| \leq 
\ltriple| A_{\alpha}(S, T)X \rtriple|. 
$$
\end{Prop}

{\it Proof:}
We firstly calculate
\[
\frac{H_{\alpha }( e^t,1 )}{G_{\alpha }( e^t,1 )} 
= \frac{2\alpha e^{\alpha t}}{e^{\alpha t} + 1}
\frac{( e^t - 1 )}{e^{\alpha t} - 1}
\frac{e^{\alpha t} - 1}{\alpha e^{\alpha t/2}( e^t - 1 )} 
= \frac{2 e^{\alpha t/2}}{e^{\alpha t} + 1} 
= \frac{2}{e^{\alpha t/2} + e^{ - \alpha t/2}}
 = \frac{1}{\cosh \frac{\alpha t}{2}},
\]
which is a positive definite function. Thus we have $H_{\alpha} \preceq G_{\alpha}$ so that the first inequality of this proposition thanks to Theorem \ref{Hiai-Kosaki-theorem}.

The calculation
\[
\frac{G_{\alpha }( e^t,1 )}{L_{\alpha }( e^t,1 )} 
= \frac{\alpha e^{\alpha t/2}( e^t - 1 )}{e^{\alpha t} - 1}
 \cdot \frac{t}{e^t - 1} 
= \frac{\alpha t}{e^{\alpha t/2} - e^{ - \alpha t/2}} 
= \frac{ \frac{\alpha t}{2}}{\sinh \frac{\alpha t}{2}}
\]
implies $G_{\alpha} \preceq L_{\alpha}$. Thus we have the second inequality of this proposition.

Finally the calculation
\[
\frac{L_{\alpha }( e^t,1 )}{A_{\alpha }( e^t,1 )} 
= \frac{e^t - 1}{t} \cdot \frac{2( e^{\alpha t} - 1 )}{\alpha (e^{\alpha t}+1)(e^t - 1 )} 
= \frac{2}{{\alpha t}} \cdot \frac{{{e^{\alpha t/2}} - {e^{ - \alpha t/2}}}}{{{e^{\alpha t/2}} + {e^{ - \alpha t/2}}}} = \frac{{\tanh \frac{{\alpha t}}{2}}}{{\frac{{\alpha t}}{2}}}\]
implies  $L_{\alpha} \preceq A_{\alpha}$. Thus we have the third inequality of this proposition.

\hfill \qed

In the papers \cite{HK1,HK2}, the unitarily invariant norm inequalities of the power difference mean (or A-L-G interpolating mean) $M_{\alpha}(x,y)$ was systematically studied.
We give the relation our means with the power difference mean:
\[{M_\alpha }\left( {x,y} \right) \equiv \left\{ \begin{array}{l}
\frac{{\alpha  - 1}}{\alpha } \cdot \frac{{{x^\alpha } - {y^\alpha }}}{{{x^{\alpha  - 1}} - {y^{\alpha  - 1}}}}, \,\,\,\,\,\, \left( {x \ne y} \right)\\
\,\,\,\,x,\,\,\,\,\,\,\,\,\,\,\,\,\,\,\,\,\,\,\,\,\,\,\,\,\,\,\,\,\,\,\,\,\,\,\,\,\,\,\,\,\left( {x = y} \right).
\end{array} \right.\]

\begin{The} \label{theorem-M-G-L-A-M-section2}
For any $S, T, X \in B(\fH)$ with $S, T \geq 0$,  $ m \in \mathbb{N}$ and any unitarily invariant norm $\ltriple| \cdot \rtriple|$, we have
$$
\hspace*{-18mm}
\ltriple| M_{\frac{m}{m+1}}(S, T)X\rtriple| \leq 
\ltriple| G_{\frac{1}{m}}(S, T)X\rtriple| \leq 
\ltriple| L(S, T)X\rtriple| \leq 
\ltriple| A_{\frac{1}{m}}(S, T)X\rtriple| \leq 
\ltriple| M_{\frac{m+1}{m}}(S, T)X\rtriple| 
$$
\end{The}

{\it Proof:}
The second inequality and the third inequality have already been proven in Proposition \ref{prop-ineq-section2}.

To prove the first inequality, for $0 < \alpha, \beta <1$ we calculate 
\begin{eqnarray*}
\frac{M_{\beta }( e^{2t},1 )}{G_{\alpha }( e^{2t},1 )}
 &=& \frac{\beta  - 1}{\beta }
 \cdot \frac{e^{2\beta t} - 1}{e^{2( \beta  - 1 )t} - 1} 
 \cdot \frac{e^{2\alpha t} - 1}{\alpha e^{\alpha t}( e^{2t} - 1 )} 
= \frac{1 - \beta }{\alpha \beta } 
\cdot \frac{\sinh \beta t}{\sinh t} 
\cdot \frac{\sinh \alpha t}{\sinh ( 1 - \beta )t}\\
&& = \frac{2( 1 - \beta  )}{\alpha \beta } 
\cdot \frac{\sinh \beta t\cosh \frac{\alpha t}{2}}{\sinh t} 
\cdot \frac{\sinh \frac{\alpha t}{2}}{\sinh ( 1 - \beta  )t}.
\end{eqnarray*}
By Proposition 5 in \cite{BK},
the function $\frac{{\sinh \beta t\cosh \frac{{\alpha t}}{2}}}{{\sinh t}}$
is positive definite, if $\beta + \frac{\alpha}{2} \leq 1$ and $\frac{\alpha}{2} \leq \frac{1}{2}$.
The function $\frac{{\sinh \frac{{\alpha t}}{2}}}{{\sinh \left( {1 - \beta } \right)t}}$ is also  positive definite, if $\frac{\alpha}{2} \leq 1- \beta$.
The case $\alpha = \frac{1}{m}$ and $\beta = \frac{m}{m+1}$ satisfies the above conditions. Thus we have $M_{\frac{m}{m+1}} \preceq G_{\frac{1}{m}}$ which leads to the first inequality of this proposition.

To prove the last inequality, for $0<\alpha <1$ and $\beta >1$, we also calculate
 \begin{eqnarray*}
\frac{A_{\alpha }( e^{2t},1 )}{M_{\beta }( e^{2t},1 )} &=& 
\frac{\alpha \beta }{2( \beta  - 1 )} 
\cdot \frac{\sinh t\sinh ( \beta  - 1 )t}{\tanh \alpha t\sinh \beta t}
 = \frac{\alpha \beta }{2( \beta  - 1 )} 
\cdot \frac{\sinh t\cosh \alpha t\sinh ( \beta  - 1 )t}{\sinh \beta t\sinh \alpha t}\\
&=&\frac{\alpha \beta }{2( \beta  - 1 )}
 \cdot \frac{\sinh \frac{1}{\beta }( \beta t )\cosh \frac{\alpha }{\beta }( {\beta t} )}{\sinh \beta t} \cdot \frac{\sinh ( \beta  - 1 )t}{\sinh \alpha t}.
\end{eqnarray*}
By Proposition 5 in \cite{BK},
the function $\frac{{\sinh \frac{1}{\beta }\left( {\beta t} \right)\cosh \frac{\alpha }{\beta }\left( {\beta t} \right)}}{{\sinh \beta t}}$ is positive definite, if
$\frac{1}{\beta} + \frac{\alpha}{\beta} \leq 1$ and $\frac{\alpha}{\beta} \leq \frac{1}{2}$ .
The function $\frac{{\sinh \left( {\beta  - 1} \right)t}}{{\sinh \alpha t}}$ is also positive definite, if $\beta -1 \leq \alpha$.
From these conditions, we have $\beta = \alpha +1$ and $\alpha \leq 1$.
The case $\alpha = \frac{1}{m}$ and $\beta =\frac{m+1}{m}$ satisfies the above conditions. Thus we have $A_{\frac{1}{m}} \preceq M_{\frac{m+1}{m}}$  which leads to the last inequality.

\hfill \qed

\begin{Rem}
Since $\frac{m+1}{m} < \frac{m}{m-1}$,  by Theorem 2.1 in \cite{HK2}, we have
 $ M_{\frac{m+1}{m}}\preceq M_{\frac{m}{m-1}}$. Thus we have 
$$\ltriple| M_{\frac{m+1}{m}}(S, T)X \rtriple| \leq \ltriple| M_{\frac{m}{m-1}}(S, T)X \rtriple|,$$
which means Theorem \ref{theorem-M-G-L-A-M-section2} gives a general result for Theorem \ref{FY_theorem2013}. At the same time, the second inequality and the third one give tighter bounds than the results given in  Theorem \ref{HK_theorem1999}.
\end{Rem}

\begin{Prop} \label{prop-H-M-section2}
For any $S, T, X \in B(\fH)$ with $S, T \geq 0$, $m =1,2$ and any unitarily invariant norm $\ltriple| \cdot \rtriple|$, we have
$$
\ltriple| H_{\frac{1}{m}}(S, T)X \rtriple| \leq
\ltriple| M_{\frac{m}{m+1}}(S, T)X \rtriple|.
$$
\end{Prop}

{\it Proof:}
For $0 < \alpha, \beta <1$ we calculate,
\begin{eqnarray*}
\frac{H_{\alpha }( e^{2t},1 )}{M_{\beta }( e^{2t},1 )} 
&=& \frac{\alpha \beta }{\beta  - 1}
\frac{\sinh t}{\cosh \alpha t\sinh \alpha t}
\frac{\sinh ( \beta  - 1 )t}{\sinh \beta t} 
= \frac{\alpha \beta }{1 - \beta }\frac{\sinh t\sinh (1 - \beta )t}{\cosh \alpha t\sinh \alpha t\sinh \beta t}\\
&=& \frac{2\alpha \beta }{1 - \beta }\frac{\sinh t\sinh ( 1 - \beta )t}{\sinh 2\alpha t\sinh \beta t}.
\end{eqnarray*}
The function  $\frac{{\sinh t}}{{\sinh 2\alpha t}}$ is positive definite, if $1 \leq 2\alpha$. The function $\frac{{\sinh \left( {1 - \beta } \right)t}}{{\sinh \beta t}}$
is also positive definite, if $1-\beta \leq \beta$. Thus $\frac{{{H_\alpha }\left( {{e^{2t}},1} \right)}}{{{M_\beta }\left( {{e^{2t}},1} \right)}}$ is positive definite, if 
$\frac{1}{2} \leq \alpha,\beta <1$. The case $\alpha = \frac{1}{m}$ and $\beta =\frac{m}{m+1}$ for $m=1,2$ satisfies the  condition $\frac{1}{2} \leq \alpha,\beta <1$ so that we have this proposition.

\hfill \qed

\begin{Rem}
We do not have the scalar inequality $H_{1/3}(t,1) \leq M_{3/4}(t,1)$ for $t>0$ in general, so that Proposition \ref{prop-H-M-section2} is not true for $m=3$.
We also do not have the scalar inequality $H_{1/3}(t,1) \geq M_{3/4}(t,1)$ for $t>0$, in general.
\end{Rem}

%%%%%%%%%%%%%%%%%%%%%%%%%%%%%%%%%%%%%%%%%%%%%%%%%%%%%%%%%%%%%%%%%%%%%%%%%%%%%%%%%%%%%%%%%%%%%%%%%%%%%%%%%%%%%%%%%%%%%%%%%%%%%%%%%%
%%%%%%%%%%%%%%%%%%%%%%%%%%%%%%%%%%%%%%%%%%%%%%%%%%%%%%%%%%%%%%%%%%%%%%%%%%%%%%%%%%%%%%%%%%%%%%%%%%%%%%%%%%%%%%%%%%%%%%%%%%%%%%%%%%
%%%%%%%%%%%%%%%%%%%%%%%%%%%%%%%%%%%%%%%%%%%%%%%%%%%%%%%%%%%%%%%%%%%%%%%%%%%%%%%%%%%%%%%%%%%%%%%%%%%%%%%%%%%%%%%%%%%%%%%%%%%%%%%%%%
%%%%%%%%%%%%%%%%%%%%%%%%%%%%%%%%%%%%%%%%%%%%%%%%%%%%%%%%%%%%%%%%%%%%%%%%%%%%%%%%%%%%%%%%%%%%%%%%%%%%%%%%%%%%%%%%%%%%%%%%%%%%%%%%%%
%%%%%%%%%%%%%%%%%%%%%%%%%%%%%%%%%%%%%%%%%%%%%%%%%%%%%%%%%%%%%%%%%%%%%%%%%%%%%%%%%%%%%%%%%%%%%%%%%%%%%%%%%%%%%%%%%%%%%%%%%%%%%%%%%%
%%%%%%%%%%%%%%%%%%%%%%%%%%%%%%%%%%%%%%%%%%%%%%%%%%%%%%%%%%%%%%%%%%%%%%%%%%%%%%%%%%%%%%%%%%%%%%%%%%%%%%%%%%%%%%%%%%%%%%%%%%%%%%%%%%
%%%%%%%%%%%%%%%%%%%%%%%%%%%%%%%%%%%%%%%%%%%%%%%%%%%%%%%%%%%%%%%%%%%%%%%%%%%%%%%%%%%%%%%%%%%%%%%%%%%%%%%%%%%%%%%%%%%%%%%%%%%%%%%%%%
%%%%%%%%%%%%%%%%%%%%%%%%%%%%%%%%%%%%%%%%%%%%%%%%%%%%%%%%%%%%%%%%%%%%%%%%%%%%%%%%%%%%%%%%%%%%%%%%%%%%%%%%%%%%%%%%%%%%%%%%%%%%%%%%%%
%%%%%%%%%%%%%%%%%%%%%%%%%%%%%%%%%%%%%%%%%%%%%%%%%%%%%%%%%%%%%%%%%%%%%%%%%%%%%%%%%%%%%%%%%%%%%%%%%%%%%%%%%%%%%%%%%%%%%%%%%%%%%%%%%%
%%%%%%%%%%%%%%%%%%%%%%%%%%%%%%%%%%%%%%%%%%%%%%%%%%%%%%%%%%%%%%%%%%%%%%%%%%%%%%%%%%%%%%%%%%%%%%%%%%%%%%%%%%%%%%%%%%%%%%%%%%%%%%%%%%
%%%%%%%%%%%%%%%%%%%%%%%%%%%%%%%%%%%%%%%%%%%%%%%%%%%%%%%%%%%%%%%%%%%%%%%%%%%%%%%%%%%%%%%%%%%%%%%%%%%%%%%%%%%%%%%%%%%%%%%%%%%%%%%%%%
%%%%%%%%%%%%%%%%%%%%%%%%%%%%%%%%%%%%%%%%%%%%%%%%%%%%%%%%%%%%%%%%%%%%%%%%%%%%%%%%%%%%%%%%%%%%%%%%%%%%%%%%%%%%%%%%%%%%%%%%%%%%%%%%%%

\section{Norm continuity in parameter}
In this section, we consider the norm continuity argument with respect to the parameter on our introduced means. 
Since we have the relation $H_{\alpha}(x,y) = G_{2\alpha}(x,y)$, we firstly consider the norm continuity in parameter on $G_{\alpha}(S,T)$.

\begin{Prop}
Let $S, T, X \in B(\fH)$ with $S, T \geq 0$. 
If $0 \leq \alpha < \beta \leq 2$ and $\ltriple| G_{\alpha}(S,T)X \rtriple| < \infty $, then we have for any unitarily invariant norm $\ltriple|\cdot \rtriple|$,
$$
\lim_{\beta \rightarrow \beta '} \ltriple| G_{\beta}(S,T)X -G_{\beta '}(S,T)X  \rtriple| =0.
$$
\end{Prop}
{\it Proof:}
From the following equality (See Eq.(1.4) in \cite{HK2} for example.)
$$
\frac{G_{\beta}(e^{2t},1)}{G_{\alpha}(e^{2t},1)} = \frac{\beta}{\alpha} \cdot \frac{\sinh \alpha t}{\sinh\beta t} =  \frac{\beta}{\alpha}  \int_{-\infty}^{\infty}
e^{its}\frac{\sin\left(\frac{\pi \alpha}{\beta}\right)}{2\beta \left\{ \cosh\left(\frac{\pi s}{\beta}\right) +\cos\left(\frac{\pi \alpha}{\beta}\right)\right\}} ds,
$$
we have for $0\leq \alpha <\beta \leq 2$, 
$$
G_{\beta}(S,T)X = \int_{-\infty}^{\infty} \left(S_{supp S}\right)^{ix}\left(G_{\alpha}(S,T)X\right)\left(T_{supp T}\right)^{-ix} \frac{\sin\left(\frac{\pi \alpha}{\beta}\right)}{2\alpha \left\{ \cosh\left(\frac{\pi s}{\beta}\right) +\cos\left(\frac{\pi \alpha}{\beta}\right)\right\}}dx,
$$
applying Theorem 3.4 in \cite{Hiai-Kosaki-book} with $G_{\beta}(1,0)=0$.
Where $S_{supp S}$ represents the support projection of $S$.
Thus we have
\begin{eqnarray*}
&& \hspace*{-16mm} \ltriple| G_{\beta}(S,T)X -G_{\beta '}(S,T)X \rtriple| \leq
\lVert \frac{\sin\left(\frac{\pi \alpha}{\beta}\right)}{2\alpha \left\{ \cosh\left(\frac{\pi s}{\beta}\right) +\cos\left(\frac{\pi \alpha}{\beta}\right)\right\}}-\frac{\sin\left(\frac{\pi \alpha}{\beta '}\right)}{2\alpha \left\{ \cosh\left(\frac{\pi s}{\beta '}\right) +\cos\left(\frac{\pi \alpha}{\beta '}\right)\right\}} \rVert_{1}\\
&& \hspace*{4cm}\times \ltriple|  G_{\alpha}(S,T)X \rtriple| \rightarrow 0
\qquad(\beta \rightarrow \beta '),
\end{eqnarray*}
by the Lebesgue dominated convergence theorem.

\hfill \qed

We secondly consider the norm continuity in parameter on $A_{\alpha}(S,T)$.
\begin{Prop}  \label{section3_prop02}
Let  $S, T, X \in B(\fH)$ with $S, T \geq 0$. 
If $0 < \alpha < \beta \leq 1$, then we have for any unitarily invariant norm $\ltriple| \cdot \rtriple|$,
\begin{equation} \label{prop01-ineq01-AM-section3}
\ltriple| A_{\alpha} (S, T) X \rtriple| \leq \ltriple| A_{\beta} (S, T) X \rtriple| \leq \frac{2\beta -\alpha}{\alpha} \ltriple| A_{\alpha} (S, T) X \rtriple|
\end{equation}
and
\begin{equation} \label{prop01-ineq02-AM-section3}
\ltriple|  A_{\alpha} (S, T)X -A_{\beta}(S, T)X \rtriple| \leq 
\frac{2(\beta -\alpha)}{\alpha}\ltriple| A_{\alpha} (S, T)X \rtriple|.
\end{equation}
\end{Prop}

{\it Proof:}
The first inequality of (\ref{prop01-ineq01-AM-section3}) has been proved in (iii) of Proposition \ref{prop-monotonicity-section2}.
Since $\frac{1}{\cosh\alpha t}$ and $\frac{\sinh(\beta -\alpha)t}{\sinh \beta t}$ are positive definite functions,
$$
1-\frac{\alpha}{\beta}\cdot \frac{A_{\beta}(e^{2t},1)}{A_{\alpha}(e^{2t},1)} =
\frac{\cosh \alpha t \sinh \beta t- \cosh \beta t \sinh \alpha t}{\cosh \alpha t \sinh \beta t} = \frac{1}{\cosh\alpha t}\cdot \frac{\sinh(\beta -\alpha)t}{\sinh \beta t}
$$
is positive definite. If we set
$$
A(s,t) \equiv \frac{\beta}{\beta -\alpha} A_{\alpha}(s,t) -\frac{\alpha}{\beta -\alpha} A_{\beta}(s,t),
$$
then we have
$
\frac{A(e^t,1)}{A_{\alpha}(e^t,1)}=\frac{\beta}{\beta -\alpha}\cdot \frac{1}{\cosh\frac{\alpha t}{2}}\cdot \frac{\sinh\frac{(\beta -\alpha)t}{2}}{\sinh\frac{\beta t}{2}},
$
which is a positive definite function. Thus we have 
\begin{equation} \label{prop01-proof01-AM-section3}
\ltriple| \frac{\beta}{\beta -\alpha} A_{\alpha}(S, T)X -\frac{\alpha}{\beta -\alpha} A_{\beta}(S, T)X \rtriple| \leq \ltriple| A_{\alpha}(S, T)X \rtriple| 
\end{equation}
by Theorem 2.4 in \cite{Kosaki}. (Actually, $A(s,t)$ may not be a symmetric homogeneous mean. However we easily find that it satisfies $A(s,t) = A(t,s)$ and $A(s,s) =s$. Then Theorem 2.4 in \cite{Kosaki} assures that the inequality (\ref{prop01-proof01-AM-section3}) is valid.)  
Therefore we have
\begin{eqnarray*}
\frac{\alpha}{\beta-\alpha}\ltriple| A_{\beta}(S, T)X \rtriple| &\leq& \frac{\beta}{\beta-\alpha}\ltriple| A_{\alpha}(S, T)X \rtriple|
+ \ltriple| \frac{\beta}{\beta -\alpha} A_{\alpha}(S, T)X -\frac{\alpha}{\beta -\alpha} A_{\beta}(S, T)X \rtriple| \\
&\leq& \left(\frac{\beta}{\beta -\alpha} +1 \right) \ltriple| A_{\alpha}(S, T)X  \rtriple| = \frac{2\beta -\alpha}{\beta -\alpha} \ltriple| A_{\alpha} (S, T) X \rtriple|,
\end{eqnarray*}
which is the second inequality of  (\ref{prop01-ineq01-AM-section3}).

We prove the inequality  (\ref{prop01-ineq02-AM-section3}).
\begin{eqnarray}
&&\ltriple| A_{\alpha} (S, T)X -A_{\beta}(S, T)X \rtriple| =
\ltriple| \left(1-\frac{\beta}{\alpha}\right) A_{\alpha}(S, T)X+\frac{\beta}{\alpha}   A_{\alpha}(S, T)X -A_{\beta}(S, T)X \rtriple| \nonumber \\
&&\hspace*{3cm}\leq \left( \frac{\beta -\alpha}{\alpha}\right) \ltriple|   A_{\alpha}(S, T)X \rtriple| +\ltriple|  \frac{\beta}{\alpha}   A_{\alpha}(S, T)X -A_{\beta}(S, T)X  \rtriple|.
\label{prop01-proof02-AM-section3}
\end{eqnarray}
From the inequality (\ref{prop01-proof01-AM-section3}), we have
$$
\ltriple| \frac{\beta}{\alpha}A_{\alpha}(S, T)X -A_{\beta}(S, T)X \rtriple| \leq \frac{\beta -\alpha}{\alpha}\ltriple| A_{\alpha}(S, T)X \rtriple|.
$$
Thus the right hand side of the inequality  (\ref{prop01-proof02-AM-section3}) is bounded from the above:
$$
 \left( \frac{\beta -\alpha}{\alpha}\right) \ltriple| A_{\alpha}(S, T)X \rtriple| +\ltriple| \frac{\beta}{\alpha}   A_{\alpha}(S, T)X -A_{\beta}(S, T)X  \rtriple| \leq \frac{2(\beta -\alpha)}{\alpha}\ltriple| A_{\alpha}(S, T)X \rtriple|.
$$
Thus we have the inequality (\ref{prop01-ineq02-AM-section3}).

\hfill \qed

We also have the following proposition.

\begin{Prop} \label{sec3_proposition3}
Let  $S, T, X \in B(\fH)$ with $S, T \geq 0$. 
If $0 \leq \alpha < \beta \leq 1$ and $\ltriple| A_{\beta}(S,T)X \rtriple| < \infty $, then we have for any unitarily invariant norm $\ltriple| \cdot \rtriple|$,
\begin{equation} \label{sec3_prop3}
\lim_{\alpha \rightarrow \alpha '} \ltriple| A_{\alpha}(S,T)X -A_{\alpha '}(S,T)X  \rtriple| =0.
\end{equation}
\end{Prop}

{\it Proof:}
We firstly prove Eq.(\ref{sec3_prop3}) for the case $0 < \alpha < \beta \leq 1$.
For $\alpha ' \in [\alpha,\beta)$, we have
$$
\ltriple| A_{\alpha}(S,T)X-A_{\alpha '}(S,T)X \rtriple| \leq
\frac{2(\alpha '-\alpha)}{\alpha}\ltriple|  A_{\alpha}(S,T)X \rtriple|  \leq
\frac{2(\alpha '-\alpha)}{\alpha}\ltriple|  A_{\beta}(S,T)X \rtriple|.
$$
by Proposition \ref{section3_prop02}. For $\alpha '\in [\frac{\alpha}{2},\alpha]$,  we similarly have 
$$
\ltriple| A_{\alpha}(S,T)X-A_{\alpha '}(S,T)X \rtriple| \leq
\frac{2(\alpha-\alpha ')}{\alpha '}\ltriple|A_{\alpha '}(S,T)X \rtriple|  \leq
\frac{4(\alpha '-\alpha)}{\alpha}\ltriple|A_{\beta}(S,T)X \rtriple|.
$$
We thus obtain for $\alpha '\in[\frac{\alpha}{2},\beta)$,
$$
\ltriple|A_{\alpha}(S,T)X-A_{\alpha '}(S,T)X \rtriple| \leq
\frac{4\vert \alpha '-\alpha\vert}{\alpha}\ltriple|  A_{\beta}(S,T)X \rtriple| 
$$
which implies Eq.(\ref{sec3_prop3}) for the case $0 < \alpha < \beta \leq 1$.

We secondly show Eq.(\ref{sec3_prop3}) for the case $\alpha = 0$.
When $0<\alpha <\beta \leq 1$, we have
$$
\frac{A_{\alpha}(e^{2t},1)}{A_{\beta}(e^{2t},1)} =\frac{\alpha}{\beta}\cdot \frac{\sinh(\beta t) \cosh(\alpha t)}{\cosh(\beta t)\sinh(\alpha t)}
=\frac{\alpha}{\beta}+\frac{\alpha}{\beta}\cdot \frac{\sinh((\beta-\alpha) t)}{\cosh(\beta t)\sinh(\alpha t)}.
$$
If we put $B(s,t) \equiv A_{\alpha}(s,t) -\frac{\alpha}{\beta} A_{\beta}(s,t)$, then we have
$$
\frac{B(e^{2t},1)}{A_{\beta}(e^{2t},1)}=\frac{\alpha}{\beta}\cdot \frac{\sinh((\beta-\alpha) t)}{\cosh(\beta t)\sinh(\alpha t)}
$$
which is a positive definite function as shown in (iii) of Proposition \ref{prop-monotonicity-section2}. We also find that
$$
\frac{A_0(e^{2t},1)}{A_{\beta}(e^{2t},1)}=\frac{1}{\beta t}\cdot \frac{\sinh(\beta t)}{\cosh(\beta t)}
$$
in the limit $\alpha \rightarrow 0$.
Then we put the Fourier transforms $\hat{\phi}_{\alpha,\beta}(t)$ and
$\hat{\phi}_{0,\beta}(t)$ of two  functions $\phi_{\alpha,\beta}(s)$ and $\phi_{0,\beta}(s)$ in the following:
\begin{eqnarray*}
&&\int_{-\infty}^{\infty} e^{ist}\phi_{\alpha,\beta}(s)ds = \hat{\phi}_{\alpha,\beta}(t) \equiv \frac{\alpha}{\beta}\cdot \frac{\sinh((\beta-\alpha) t)}{\cosh(\beta t)\sinh(\alpha t)},\\
&&\int_{-\infty}^{\infty} e^{ist}\phi_{0,\beta}(s)ds = \hat{\phi}_{0,\beta}(t) \equiv 
\frac{1}{\beta t}\cdot \frac{\sinh(\beta t)}{\cosh(\beta t)}.
\end{eqnarray*}
Since we have $B(1,0)=0$ and $A_{0}(1,0) =0$, we have
\begin{eqnarray*}
&& A_{\alpha}(S,T)X-\frac{\alpha}{\beta} A_{\beta}(S,T)X = \int_{-\infty}^{\infty}(S_{suppS})^{is}(A_{\beta}(S,T)X)(T_{suppT})^{-is} \phi_{\alpha,\beta}(s)ds\\
&& A_{0}(S,T)X= \int_{-\infty}^{\infty}(S_{suppS})^{is}(A_{\beta}(S,T)X)(T_{suppT})^{-is} \phi_{0,\beta}(s)ds
\end{eqnarray*}
from Theorem 3.4 in \cite{Hiai-Kosaki-book}.
Then we have
%\begin{eqnarray*}
%&&A_{\alpha}(S,T)X-A_{0}(S,T)X = \frac{\alpha}{\beta} A_{\beta}(S,T)X \\
%&& \hspace*{3cm}+\int_{-\infty}^{\infty}(S_{suppS})^{is}(A_{\beta}(S,T)X)%(T_{suppT})^{-is} (\phi_{\alpha,\beta}(s)-\phi_{0,\beta}(s))ds.
%\end{eqnarray*}
%which implies
$$
\ltriple| A_{\alpha}(S,T)X-A_{0}(S,T)X \rtriple| \leq
\left( \frac{\alpha}{\beta} + \lVert   \phi_{\alpha,\beta} -\phi_{0,\beta} \rVert_1\right)
\ltriple|  A_{\beta}(S,T)X \rtriple|. 
$$
To prove $\lim_{\alpha \rightarrow 0+} \lVert \phi_{\alpha,\beta} -\phi_{0,\beta} \rVert_1 =0$, we have only to prove $\lim_{\alpha \rightarrow 0+} \lVert \hat{\phi}_{\alpha,\beta} -\hat{\phi}_{0,\beta} \rVert_2 =0$ thanks to Lemma 5.8 in \cite{Hiai-Kosaki-book}.
Since we have $\int_{-\infty}^{\infty} \phi_{\alpha,\beta}(s) ds = \hat{\phi}_{\alpha,\beta}(0) = \frac{\beta -\alpha}{\beta}$,
we have
$\int_{-\infty}^{\infty} \phi_{0,\beta}(s) ds = \hat{\phi}_{0,\beta}(0) = 1$
in the limit $\alpha \rightarrow 0$.
From the fact $\sinh x \geq x $ for $x \geq 0$, we also have 
$0 \leq \hat{\phi}_{\alpha,\beta}(t), \hat{\phi}_{0,\beta}(t) \leq \frac{1}{\beta | t|}$. Since $\hat{\phi}_{\alpha,\beta}(t)$ and $\hat{\phi}_{0,\beta}(t)$ are
positive definite functions, we have $\hat{\phi}_{\alpha,\beta}(t) \leq \hat{\phi}_{\alpha,\beta}(0) =\frac{\beta -\alpha}{\beta} \leq 1$
and $\hat{\phi}_{0,\beta}(t) \leq \hat{\phi}_{0,\beta}(0) =1$. (See Chapter 5 in \cite{Bhatia} for basic properties of the positive definite function.)
We thus obtain $\hat{\phi}_{\alpha,\beta}(t), \hat{\phi}_{0,\beta}(t)  \leq \min\left(1,\frac{1}{\beta | t|}\right)$
for two $L^2$-functions $\hat{\phi}_{\alpha,\beta}$ and $\hat{\phi}_{0,\beta}$.
We finally obtain $\vert \hat{\phi}_{\alpha,\beta}(t)- \hat{\phi}_{0,\beta}(t)\vert^2 \leq 4\min(1,\frac{1}{\beta t^2})$.
Since $\min(1,\frac{1}{\beta t^2})$ is integrable and $\lim_{\alpha\rightarrow 0}\hat{\phi}_{\alpha,\beta}(t)= \hat{\phi}_{0,\beta}(t)$,
we obtain  $\lim_{\alpha \rightarrow 0+} \lVert \hat{\phi}_{\alpha,\beta} -\hat{\phi}_{0,\beta} \rVert_2 =0$ by the Lebesgue dominated convergence theorem.

\hfill \qed

We note that the assumption $\ltriple| A_{\beta}(S,T)X \rtriple| < \infty$ for some $\beta \in (0,1]$ is equivalent to $\ltriple| SX+XT \rtriple| <\infty$, since we have
$\ltriple| A_{\beta}(S,T)X \rtriple| \leq \ltriple| A_{1}(S,T)X \rtriple| \leq \frac{2-\beta}{\beta}\ltriple| A_{\beta}(S,T)X \rtriple| $
using the inequality (\ref{prop01-ineq01-AM-section3}).

\section{Conclusion}
We obtained new and tight bounds for the logarithmic mean for untarily invariant norm. Our results improved the famous inequalities by Hiai and Kosaki \cite{HK1,HK2}. Concluding this paper, we summarize Theorem \ref{theorem-M-G-L-A-M-section2} by the familiar form. From the calculations
$$G_{1/m_1}(s,t)=\frac{1}{m_1} \sum_{k=1}^{m_1} s^{(2k-1)/2m_1} t^{(2m_1-(2k-1))/2m_1}$$ and 
$$A_{1/m_2}(s,t)=\frac{1}{m_2}\left( \sum_{k=0}^{m_2} s^{k/m_2} t^{(m_2-k)/m_2} -\frac{1}{2}(s+t)\right),$$ we have
$$
G_{1/m_1} (S,T)X =\frac{1}{m_1} \sum_{k=1}^{m_1} S^{(2k-1)/2m_1} XT^{(2m_1-(2k-1))/2m_1}
$$
and
$$
A_{1/m_2} (S,T)X = \frac{1}{m_2} \left(\sum_{k=0}^{m_2} S^{k/m_2}XT^{(m_2-k)/m_2} -\frac{1}{2}(SX+XT)\right).
$$
In addition, from the paper \cite{HK2}, we know that
$$
M_{m_1/(m_1+1)}(S,T)X =\frac{1}{m_1} \sum_{k=1}^{m_1}S^{k/(m_1+1)}XT^{(m_1+1-k)/(m_1+1)}
$$
and
$$
M_{m_2/(m_2-1)} (S,T)X =\frac{1}{m_2} \sum_{k=0}^{m_2-1} S^{k/(m_2-1)}XT^{(m_2-1-k)/(m_2-1)}.
$$
Thus Theorem \ref{theorem-M-G-L-A-M-section2} can be rewitten as the following inequalities which are our main result of the present paper. 
\begin{eqnarray}
\hspace*{-12mm}\frac{1}{m_1}\ltriple| \sum_{k=1}^{m_1}S^{k/(m_1+1)}XT^{(m_1+1-k)/(m_1+1)}\rtriple| &\leq&\frac{1}{m_1} \ltriple| \sum_{k=1}^{m_1} S^{(2k-1)/2m_1} XT^{(2m_1-(2k-1))/2m_1} \rtriple| \nonumber \\
&\leq& \ltriple| \int_0^1 S^{\nu}XT^{1-\nu}d\nu \rtriple| \nonumber \\
&\leq& \frac{1}{m_2} \ltriple| \sum_{k=0}^m S^{k/m_2}XT^{(m_2-k)/m_2} -\frac{1}{2}(SX+XT)\rtriple| \nonumber \\
&\leq& \frac{1}{m_2} \ltriple| \sum_{k=0}^{m_2-1} S^{k/(m_2-1)}XT^{(m_2-1-k)/(m_2-1)}\rtriple|,
\end{eqnarray}
for $S,T,X \in B(\fH)$ with $S,T \geq 0$, $m_1 \geq 1$, $m_2 \geq 2$, and any unitarily invariant norm $\ltriple| \cdot \rtriple|$.

We have also shown some properties for our means such as monotonicities and norm continuities in parameter. 

\section*{Competing interests}
The authors declare that they have no competing interests.

\section*{Acknowledgements}
%The author would like to express his deepest gratitude to Professor Fumio Hiai %for letting the author know the powerful method (Theorem \ref{Hiai-Kosaki-theorem}) and its use in great detail. The continuity result for $\alpha=0$ in %Propsition \ref{sec3_proposition3} is due to Professor Hideki Kosaki %\cite{Kosaki2013_kyoto}, and the author also would like to express his deepest %gratitude to Professor Hideki Kosaki for giving the author valuable comments %throughout the paper (especially on section 3).
I would like to  express my deepest gratitude to Professor Fumio Hiai and Professor Hideki Kosaki for giving me valuable comments to improve this manuscript.
I also was partially supported by JSPS KAKENHI Grant Number 24540146.

\end{document}